\documentclass[11pt]{article}
\usepackage{amsmath}
\usepackage{amssymb}
\usepackage{gastex}
\newtheorem{lemma}{{\sc Lemma}}[section]
\newtheorem{theo}{{\sc Theorem}}

\newtheorem{cor}{{\sc Corollary}}
\newtheorem{defn}{{\em Definition}}
\newtheorem{rem}{{\em Remark}}

\newtheorem{algo}{{\sc Algorithm}}
\newtheorem{example}{{\sc Example}}

\begin{document}

\title{Contractible Hamiltonian Cycles in Polyhedral Maps}
\author{{Dipendu Maity  and Ashish Kumar Upadhyay}\\[2mm]
{\normalsize Department of Mathematics}\\{\normalsize Indian Institute of Technology Patna }\\{\normalsize Patliputra Colony, Patna 800\,013,  India.}\\
{\small \{dipendumaity, upadhyay\}@iitp.ac.in}}
\maketitle

\vspace{-5mm}

\hrule

\begin{abstract} We present a necessary and sufficient condition for existence of a contractible Hamiltonian Cycle in the edge graph of equivelar maps on surfaces. We also present an algorithm to construct such cycles. This is further generalized and shown to hold for more general maps.
\end{abstract}

{\small

{\bf AMS classification\,:} 57Q15, 57M20, 57N05.

{\bf Keywords\,:} Contractible Hamiltonian cycles, Proper Trees in Maps,
Equivelar Maps.
}

\bigskip

\hrule

\section{Introduction and Definitions}

Recall that a graph $G := (V, E)$ is a simple graph with vertex set $V$ and edge set $E$. A {\em surface} $S$ is a connected, compact, 2-dimensional manifold without boundary. A {\em map}  on a surface $S$ is an embedding of a finite graph $G$ such that the closure of components of $S \setminus G$ are $p-gonal$ 2-disks where ($p \geq 3$). The components are also called $facets$. The map $M$ is called a {\em polyhedral map} if  non - empty intersection of any two facets of the map is either a vertex or an edge see \cite{brehm_schulte}. We call $G$ the edge graph of the map and denote it by $EG(M)$. The vertices and edges of $G$ are also called vertices and edges of the map, respectively. In what follows we will use the terms map and polyhedral map interchangeably to denote a polyhedral map. A polyhedral map $M$ is called {\em $\{p, q\}$ equivelar}, $p, q \geq 3$ if each vertex in $M$ is incident with exactly $q$ numbers of $p$-gonal facets. If $p = 3$ then the map is called a $q$ - equivelar triangulation or a degree - regular triangulation of type $q$. A path $P$ in a graph $G$ is a subgraph $P:[v_1 v_2 \ldots v_n]$ of $G$, such that the vertex set of $P$ is $V(P) = \{v_1, v_2, \ldots, v_n\}\subseteq V(G)$ and $v_{i}v_{i + 1}$ are edges in $P$ for $1 \leq i \leq n - 1$. A path $P:[v_1, v_2, \ldots, v_n]$ in $G$ is said to be a cycle if $v_{n}v_1$ is also an edge in $P$. A graph without any cycles is called a tree. Length $l(P)$ of path $P$ is the number of edges in $P$. See \cite{mohar thomassen} for details about graphs on surfaces and \cite{BondyMurthy} for graph theory related terminology.

In this article we are interested in finding out whether a Hamiltonian cycle exists in the edge graph of a polyhedral map?. Such cycles in planar graphs have been extensively studied. For example, in \cite{Tutte} Tutte showed that every 4-connected planar graph has a Hamiltonian cycle. In 1970, Gr$\ddot{u}$nbaum, \cite{grunbaum} conjectured that every 4-connected graph which admits an embedding in the torus has a Hamiltonian cycle. In \cite{fuji_naka_ozeki} it is shown that a 3-connected bipartite graph embeddable in torus has a Hamiltonian cycle if it is balanced and each vertex of one of its partite sets has degree four. We are thus led to think about vertex degree considerations in graphs while looking for such cycles. The work in this direction has been going on for quite some time. A. Altshuler \cite{a1}, \cite{a2} studied Hamiltonian cycles in the edge graph of equivelar maps on the torus. He showed that in the graph consisting of vertices and edges of equivelar maps of above type there exists a Hamiltonian cycle. Continuous analogues of such cycles in maps may or may not be homotopic to the generators of fundamental group of the surface on which they lie. The cycles which are homotopic to a generator are essential cycles and those which are homotopic to a point are inessential cycles. We are trying to figure out such inessential cycles combinatorially.

We will call a cycle in the edge graph of a map to be {\em contractible} if it bounds a $2$-disk (2-cell). For example boundary cycle of a facet is contractible. Study of such cycles are being actively pursued and they are being applied to obtain important results. For example, in \cite{mori_nakamoto} the authors use contractible Hamiltonian cycle in triangulations of projective plane to determine a linear upper bound on the number of diagonal flips to mutually transform any two triangulations into one another. In \cite{brunet_nakamoto_negami} the authors produce a contractible Hamiltonian cycle in every 5-connected triangulation of the Klein bottle.   In \cite{Upadhyay} the second author presented a necessary and sufficient condition for existence of contractible Hamiltonian cycles in equivelar triangulation of surfaces. In this article we extend the results of \cite{Upadhyay} to present a necessary and sufficient condition for existence of a contractible Hamiltonian cycle in edge graph of an equivelar map. In section \ref{algorithm} we give computer implementable steps of the ideas used in proof to determine (possibly!!) all contractible Hamiltonian cycles in a given map. This is extended in section \ref{pmap} to maps other than the equivelar maps.

We begin with some definitions which will be needed in the course of proof of main Theorem \ref{thm1}. For more details on these topics one may refer to \cite{mcmullen schulte}.

If $v$ is a vertex of a map $K$ then the number of edges incident with $v$ is called the {\em degree} of $v$ and is denoted by $\deg(v)$. If the number of vertices, edges and facets of $K$ are denoted by $f_0(K), f_1(K)$ and $f_2(K)$ respectively, then the integer $\chi(K) = f_0(K) - f_1(K) + f_2(K)$ is called the {\em Euler characteristic of $K$}. The {\em dual map} $M$, of $K$ is defined to be the map on same surface as $K$ which has for its vertices the set of facets of $K$ and two vertices $u_1$ and $u_2$ of $M$ are ends of an edge of $M$ if the corresponding facets in $K$ have an edge in common. The well known maps of type $\{3, 6\}$ and $\{6, 3\}$ on the surface of torus are examples of mutually dual maps. Let $K'$ be a subset of set of faces of $K$ and $D = \displaystyle\cup_{\small{\sigma \in K'}}\sigma$. If $D$ is connected then we will call it a {\em disc} in $K$. An edge $e$ of a facet $\sigma$ in $D$ is said to be a {\em free} edge of $\sigma$, if $e$ is not contained in any other facet in $D$. The process of deleting a facet which has a free edge in a $D$ is called an {\em elementary collapse} on $D$. Applying a sequence of elementary collapses to $D$ results into another disc $\mathcal{K}$ of $K$. We say that $D$ collapses to $\mathcal{K}$. If $D$ collapses to a point then we say that $D$ is a collapsible map. It is a fact that collapsible maps are contractible [compare \cite{hachimori}, pp. 32].

Consider a $\{p, q\}$ equivelar map $K$ on a surface $S$ that has $n$ vertices. \,:

\begin{defn}\label{defn1} Let $M$ denote the dual map of $K$. Let $T := (V, E)$ denote a tree  in the edge graph $EG(M)$ of $M$. We say that $T$ is a proper tree if the following conditions hold\,:

\begin{enumerate}
\item $\displaystyle\sum_{i = 1}^{ k }\deg(v_{i}) = n + 2(k-1)$, where $V = \{v_{1},v_{2},...,v_{k}\}$, $\deg(v)$ denotes degree of $v$ in $EG(M)$
\item whenever two vertices $u_{1}$ and $u_{2}$ of $T$ lie on a face $F$ in $M$, a path $P[u_{1},u_{2}]$ joining $u_{1}$ and $u_{2}$ in the boundary $\partial{F}$ of $F$ is a part of $T$, and
\item any path $P$ in $T$ which lies in a face $F$ of $M$ is of length at most $q - 2$, where $q =$ length of $(\partial{F})$.
\end{enumerate}
\end{defn}


\begin{rem} If the map $K$ is $\{p, q\}$ equivelar then $k = \displaystyle\frac{n - 2}{p - 2}$. Thus, for an equivelar triangulation on $n$ vertices the proper tree has exactly $n - 2$ vertices.
\end{rem}

\begin{rem} Note that the disc $D$ in $K$ which is corresponding dual of the proper tree $T$ in $M$ is collapsible and therefore it is a topological $2$-disc.

\end{rem}

\begin{defn}\label{defn2} A proper tree $T$ is called an {\em admissible proper tree} if the boundary of corresponding dual $2$-disc $D$ in $K$ is a Hamiltonian Cycle in $EG(K)$.
\end{defn}

Main result of this article is\,:
\begin{theo}\label{thm1} The edge graph $EG(K)$ of an equivelar map $K$ has a contractible Hamiltonian cycle if and only if the edge graph of corresponding dual map of $K$ has a proper tree.
\end{theo}

More generally, we prove\,:
\begin{theo}\label{thm2} The edge graph $EG(K)$ of a  map $K$ on a surface has a contractible Hamiltonian cycle if and only if the edge graph of corresponding dual map of $K$ has a proper tree.
\end{theo}

These results rely on \,:
\begin{lemma}\label{lemma1} Let $M$ denote the dual map of an $n$ vertex $\{p, q\}$-equivelar map $K$ on a surface $S$. If $\frac{n - 2}{p - 2} = m$ is an integer then $M$ has an admissible proper tree on $m$ vertices.
\end{lemma}

In the next section we give an example of a dual map of a triangulation. These are $\{3, 6\}$ and $\{6, 3\}$ maps on torus. The second example is of self dual map of type $\{4, 4\}$ on the torus. There does not exist any Hamiltonian cycle in the second example. In the section following it we present some facts and properties of a proper tree and proceed to prove the main result of this article. In section \ref{algorithm} we also give computer implementable steps to find out a Hamiltonian cycle.

\section{Example and Results}
\hrule
\smallskip
\begin{example}\label{example1}
{\small {\bf $\{3, 6\}$ and $\{6, 3\}$-equivelar maps (left) and $\{4, 4\}$-equivelar map (right) on the torus }}
\end{example}
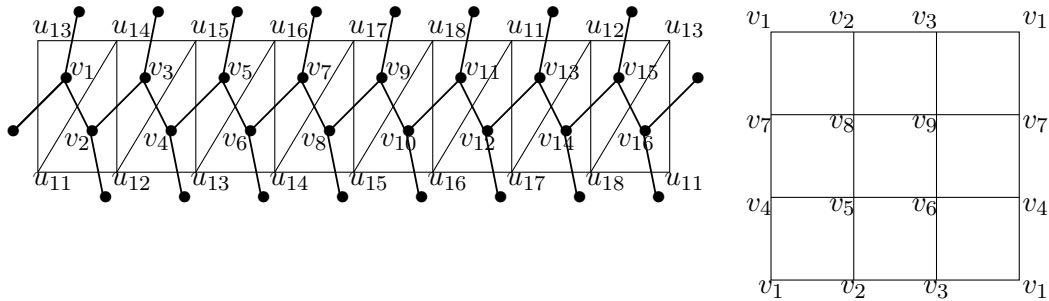
\begin{picture}(-20,55)(13,0)
\setlength{\unitlength}{3.5mm}
\drawpolygon(5,10)(29,10)(29,15)(5,15)
\put(4.8,9.5){$\scriptsize{u_{11}}$}\put(4.8,15.3){$\scriptsize{u_{13}}$}
\drawline[AHnb=0](8,10)(8, 15)
\put(7.8,9.5){$\scriptsize{u_{12}}$}\put(7.8,15.3){$\scriptsize{u_{14}}$}
\drawline[AHnb=0](11,10)(11, 15)
\put(10.8,9.5){$\scriptsize{u_{13}}$}\put(10.8,15.3){$\scriptsize{u_{15}}$}
\drawline[AHnb=0](14,10)(14, 15)
\put(13.8,9.5){$\scriptsize{u_{14}}$}\put(13.8,15.3){$\scriptsize{u_{16}}$}
\drawline[AHnb=0](17,10)(17, 15)
\put(16.8,9.5){$\scriptsize{u_{15}}$}\put(16.8,15.3){$\scriptsize{u_{17}}$}
\drawline[AHnb=0](20,10)(20, 15)
\put(19.8,9.5){$\scriptsize{u_{16}}$}\put(19.8,15.3){$\scriptsize{u_{18}}$}
\drawline[AHnb=0](23,10)(23, 15)
\put(22.8,9.5){$\scriptsize{u_{17}}$}\put(22.8,15.3){$\scriptsize{u_{11}}$}
\drawline[AHnb=0](26,10)(26, 15)
\put(25.8,9.5){$\scriptsize{u_{18}}$}
\put(25.8,15.3){$\scriptsize{u_{12}}$}
\put(28.8,9.5){$\scriptsize{u_{11}}$}
\put(28.8,15.3){$\scriptsize{u_{13}}$}
\drawline[AHnb=0](5,10)(8, 15)
\drawline[AHnb=0](8,10)(11, 15)
\drawline[AHnb=0](11,10)(14, 15)
\drawline[AHnb=0](14,10)(17, 15)
\drawline[AHnb=0](17,10)(20, 15)
\drawline[AHnb=0](20,10)(23, 15)
\drawline[AHnb=0](23,10)(26, 15)
\drawline[AHnb=0](26,10)(29, 15)
\drawline[AHnb=0,linewidth=.07](4,11.5)(6,13.5)\put(3.8,11.3){$\bullet$}
\put(6.2,13.7){$\scriptsize{v_{1}}$}
\drawline[AHnb=0,linewidth=.07](6,13.5)(7,11.5)\put(5.8,13.3){$\bullet$}
\drawline[AHnb=0,linewidth=.07](7,11.5)(9,13.5)\put(6.8,11.3){$\bullet$}
\put(6,11){$\scriptsize{v_{2}}$}\put(9.2,13.7){$\scriptsize{v_{3}}$}
\drawline[AHnb=0,linewidth=.07](9,13.5)(10,11.5)\put(8.8,13.3){$\bullet$}
\drawline[AHnb=0,linewidth=.07](10,11.5)(12,13.5)\put(9.8,11.3){$\bullet$}
\put(9,11){$\scriptsize{v_{4}}$} \put(12.2,13.7){$\scriptsize{v_{5}}$}
\drawline[AHnb=0,linewidth=.07](12,13.5)(13,11.5)\put(11.8,13.3){$\bullet$}
\drawline[AHnb=0,linewidth=.07](13,11.5)(15,13.5)\put(12.8,11.3){$\bullet$}
\put(12,11){$\scriptsize{v_{6}}$}\put(15.2,13.7){$\scriptsize{v_{7}}$}
\drawline[AHnb=0,linewidth=.07](15,13.5)(16,11.5)\put(14.8,13.3){$\bullet$}
\drawline[AHnb=0,linewidth=.07](16,11.5)(18,13.5)\put(15.8,11.3){$\bullet$}
\put(15,11){$\scriptsize{v_{8}}$}\put(18.2,13.7){$\scriptsize{v_{9}}$}
\drawline[AHnb=0,linewidth=.07](18,13.5)(19,11.5)\put(17.8,13.3){$\bullet$}
\drawline[AHnb=0,linewidth=.07](19,11.5)(21,13.5)\put(18.8,11.3){$\bullet$}
\put(18,11){$\scriptsize{v_{10}}$}\put(21.2,13.7){$\scriptsize{v_{11}}$}
\drawline[AHnb=0,linewidth=.07](21,13.5)(22,11.5)\put(20.8,13.3){$\bullet$}
\drawline[AHnb=0,linewidth=.07](22,11.5)(24,13.5)\put(21.8,11.3){$\bullet$}
\put(21,11){$\scriptsize{v_{12}}$}\put(24.2,13.7){$\scriptsize{v_{13}}$}
\drawline[AHnb=0,linewidth=.07](24,13.5)(25,11.5)\put(23.8,13.3){$\bullet$}
\drawline[AHnb=0,linewidth=.07](25,11.5)(27,13.5)\put(24.8,11.3){$\bullet$}
\put(24,11){$\scriptsize{v_{14}}$}\put(27.2,13.7){$\scriptsize{v_{15}}$}
\drawline[AHnb=0,linewidth=.07](27,13.5)(28,11.5)\put(26.8,13.3){$\bullet$}
\drawline[AHnb=0,linewidth=.07](28,11.5)(30,13.5)\put(27.8,11.3){$\bullet$}
\put(27,11){$\scriptsize{v_{16}}$}
\put(29.8,13.3){$\bullet$}
\drawline[AHnb=0,linewidth=.07](4,11.5)(6,13.5)
\drawline[AHnb=0,linewidth=.07](6,13.5)(6.5,16)\put(6.3,15.8){$\bullet$}
\drawline[AHnb=0,linewidth=.07](7,11.5)(7.5,9)\put(7.3,8.8){$\bullet$}
\drawline[AHnb=0,linewidth=.07](9,13.5)(9.5,16)\put(9.3,15.8){$\bullet$}
\drawline[AHnb=0,linewidth=.07](10,11.5)(10.5,9)\put(10.3,8.8){$\bullet$}
\drawline[AHnb=0,linewidth=.07](12,13.5)(12.5,16)\put(12.3,15.8){$\bullet$}
\drawline[AHnb=0,linewidth=.07](13,11.5)(13.5,9)\put(13.3,8.8){$\bullet$}
\drawline[AHnb=0,linewidth=.07](15,13.5)(15.5,16)\put(15.3,15.8){$\bullet$}
\drawline[AHnb=0,linewidth=.07](16,11.5)(16.5,9)\put(16.3,8.8){$\bullet$}
\drawline[AHnb=0,linewidth=.07](18,13.5)(18.5,16)\put(18.3,15.8){$\bullet$}
\drawline[AHnb=0,linewidth=.07](19,11.5)(19.5,9)\put(19.3,8.8){$\bullet$}
\drawline[AHnb=0,linewidth=.07](21,13.5)(21.5,16)\put(21.3,15.8){$\bullet$}
\drawline[AHnb=0,linewidth=.07](22,11.5)(22.5,9)\put(22.3,8.8){$\bullet$}
\drawline[AHnb=0,linewidth=.07](24,13.5)(24.5,16)\put(24.3,15.8){$\bullet$}
\drawline[AHnb=0,linewidth=.07](25,11.5)(25.5,9)\put(25.3,8.8){$\bullet$}
\drawline[AHnb=0,linewidth=.07](27,13.5)(27.5,16)\put(27.3,15.8){$\bullet$}
\drawline[AHnb=0,linewidth=.07](28,11.5)(28.5,9)\put(28.3,8.8){$\bullet$}
\end{picture}

\begin{picture}(30,12)(-85,-22)
\setlength{\unitlength}{2.2mm}
\drawpolygon(5,5)(20,5)(20,20)(5,20)
\drawline[AHnb=0](10,5)(10,20)
\drawline[AHnb=0](15,5)(15,20)
\drawline[AHnb=0](5,10)(20,10)
\drawline[AHnb=0](5,15)(20,15)
\put(4.2,4.2){$\scriptsize{v_{1}}$}
\put(9.2,4.2){$\scriptsize{v_{2}}$}
\put(14.2,4.2){$\scriptsize{v_{3}}$}
\put(20.2,4.2){$\scriptsize{v_{1}}$}
\put(3.5,9.2){$\scriptsize{v_{4}}$}
\put(3.5,14.2){$\scriptsize{v_{7}}$}
\put(3.5,20.5){$\scriptsize{v_{1}}$}
\put(8.5,9.2){$\scriptsize{v_{5}}$}
\put(13.5,9.2){$\scriptsize{v_{6}}$}
\put(20.2,9.2){$\scriptsize{v_{4}}$}
\put(8.5,14.2){$\scriptsize{v_{8}}$}
\put(13.5,14.2){$\scriptsize{v_{9}}$}
\put(20.2,14.2){$\scriptsize{v_{7}}$}
\put(8.5,20.5){$\scriptsize{v_{2}}$}
\put(13.5,20.5){$\scriptsize{v_{3}}$}
\put(20.2,20.5){$\scriptsize{v_{1}}$}
\end{picture}
\vspace{-1in}
\hrule
\section{Proper tree in equivelar maps\,:}
\begin{lemma}\label{lem1} Let $T$ be a proper tree in a $\{q, p\}$ equivelar map $M$. Then $T\bigcap F\neq \emptyset$ for any face $F$ of $M$.
\end{lemma}

\noindent{\sc Proof of Lemma}\ref{lem1} Let $V(T)$ and $E(T)$ respectively be the set of vertices and edges of the tree $T$. We construct two sets $E$ and $\widetilde{F}$ as follows. Let $E$ be a singleton set which contains a vertex $v_{1}\in V(T)$ at $1^{st}$(initial) step and $F_{1,1},F_{1,2},.....,F_{1,p}$ are the faces of $M$ such that all the faces $F_{1,1},F_{1,2},.....,F_{1,p}$ are incident at vertex $v_{1}$. Put $\widetilde{F}$ = $\{$F$_{1,1}$,F$_{1,2}$,.....,F$_{1,p}$$\}$. At $i^{th}$ step, we choose a vertex $v_{i}\in V(T)\setminus E$ such that $\{w,v_{i}\}\in E(T)$, where $w\in E$ and put in set $E$. Since $v_{i}$ be the new vertex of $E$, so there are some adjacent faces of $v_{i}$ which are not in the set $\widetilde{F}$. Now claim is, there are exactly $r = (p-2)$ faces $F_{i,1},F_{i,2},\dots,F_{i,p-2}$ where each are different from all the faces of $\widetilde{F}$ and all are incident at vertex $v$. Suppose $r\neq(p-2)$. Then there are two possibilities:

One: Suppose $r < p-2$, then there are at least three faces $F^{i}_{1},F^{i}_{2}$ and $F^{i}_{3}$ which are in $\widetilde{F}$ and all are incident at vertex $v_{i}$. Suppose $\{w,v_{i}\}$ be an common edge of $F^{i}_{1}~and~ F^{i}_{2}$ and there exist a vertex $u\in E$ such that $u\in V(F^{i}_{3})$. Since $v_{1}$ and $u$ are the vertices of tree $T$, so there exist a path $P_{1}(v_{1}\rightarrow u)$. Similarly there exist a path from $P_{2}(v_{1}\rightarrow v_{i})$. Also,  since $v_{i}$ and $u$ lie on the face $F^{i}_{3}$ and $u,v_{i}\in V(T)$. Therefore the path $P_3(v_{i}\rightarrow u)$ in $F^{i}_{3}$ also party of tree. Hence $P_{1}\bigcup  P_{2} \bigcup P_{3}$ contains a $cycle$ or collection of cycles. This can not happen because $T$ is a tree. Hence $r\geq(p-2)$.

Two: Suppose $r>(p-2)$. We know at each vertex there are exactly $p$ faces adjacent. Suppose we are choosing an edge $\{u_{1},u_{2}\}$ in $E(T)$ such that $u_{1}\in E$ and $u_{2}\in V(T)\backslash E$ and this edge is the common edge between exactly two face and both are incident at the vertices $u_{1}$, $u_{2}$. So both the faces are repetition at the vertex $u_{2}$. Hence $r\leq (p-2)$.\\Hence $r$ = $p-2$ and at each step there are $(p-2)$ faces which are different from all the faces of $\widetilde{F}$ and collect them in $\widetilde{F}$. Thus the number of faces in $\widetilde{F}=p+\underbrace{(p-2)+(p-2)+....+(p-2)}$ (repetitions of $(p-2) = (\#|V|-1)  )$. That is $\#\widetilde{F} = p+(p-2)(\frac{n-2}{p-2}-1)= p + n - p = n$. So after $\#|V|^{th}$  step $\widetilde{F}$ will contain all the faces of polyhedral map. Hence tree touches all the faces of $M$. This proves the Lemma

\smallskip

\begin{lemma}\label{lem2} Let $K$ be a $n$ vertex $\{p,q\}$ equivelar map of a surface $S$. Let $M$ denote the dual polyhedron corresponding to $K$ and $T$ be a $\frac{n-2}{p-2}$ vertex proper tree in $M$. Let $D$ denote the $subcomplex$ of $K$ which is $dual$ of $T$. Then $D$ is a $2-disk$ and $bd(D)$ is a Hamiltonian cycle in $K$.
\end{lemma}

\noindent{\sc Proof of Lemma } \ref{lem2}: By definition of a dual, $D$ consists of $\frac{n-2}{p-2}$ $p-gons$ corresponding to $\frac{n-2}{p-2}$ vertices of $T$. Two $p-gons$ in $D$ have an edge in common if the corresponding vertices are adjacent in $T$.Here the set $D$ is a collapsible $simplicial$ complex and hence it is a  $2-disk$.Since $T$ has vertices of degree one, $bd(D)\neq\emptyset$, and being boundary complex of a $2-disk$, it is a connected cycle. Observe that the number of edges in $\frac{n-2}{p-2}$ $p-gons$ is $p(\frac{n-2}{p-2}$) and for each edge of $T$ exactly $2 edges$ are identified. Hence the number of edges which remains unidentified in $D$ is $p(\frac{n-2}{p-2})-2(\frac{n-2}{p-2}-1)=\frac{pn-2p-2n+2p}{p-2}=\frac{n(p-2)}{p-2}=n$. Hence the number of vertices in $bd(D):=\partial D = n$. If the vertices $v_{1},v_{2}\in\partial D$ such that $v_{1},v_{2}$ lie on a path of length $< n$ and $v_{1}=v_{2}$. This means the faces $F_{1}$ and $F_{2}$ in $D$ with $v_{1}\in F_{1},v_{2}\in F_{2},F_{1}\neq F_{2}$ and $F_{1}$ not adjacent to $F_{2}$. Thus there exist a face $F'$ in $D$ such that the vertex $u_{F'}$ in $T$ corresponding to $F'$ does not belong to the face $F(v_{1})$ corresponding to vertex $v_{1}$. But this contradicts that $T$ is a proper tree. Thus a the cycle $\partial D$ contains exactly $n$ distinct vertices. Since $\#V(K)=n$, $\partial D$ is a Hamiltonian cycle in $K$. This proves the Lemma.

\smallskip

\begin{lemma}\label{lemm1} Let $M$ denote the dual map of an $n$ vertex$\{p, q\}$ equivelar map $K$ on a surface $S$. If $\frac{n - 2}{p - 2}$ is not an integer then $M$ does not have any admissible proper tree.
\end{lemma}

\noindent {\sc Proof of lemma} \ref{lemm1}: Let $T$ denote a tree on $m$ vertices and $V(T)=\{v_{1},v_{2},....,v_{m}\}$ be the set of vertices of  $T$. Let $\frac{n-2}{p-2}$ be as in statement of Lemma \ref{lemm1}. Then, there are two possibilities, namely, $m\leq\lfloor\frac{n-2}{p-2}\rfloor$ or $m\geq\lceil\frac{n-2}{p-2}\rceil$. As shown in the proof of Lemma \ref{lem1}, at a vertex $v_{1}$ exactly $p$ faces of $M$ are incident and if $v_i$ is adjacent to $v_1$ then exactly $p-2$ faces distinct from the $p$ faces containing $v_1$ are incident at $v_{i}(i\neq1)$. We take union of all the faces and denote this union by $\widetilde{F}$.

When $m\leq\lfloor\frac{n-2}{p-2}\rfloor$, the number of faces in $\widetilde{F} = p + \underbrace{(p-2)+(p-2)+....+(p-2)}$ ($m -1$ repetitions of $(p-2)$)
\begin{eqnarray*}
& = & p+(m-1)(p-2)\\
& \leq & p+(\lfloor\frac{n-2}{p-2}\rfloor -1)(p-2)\\
& = &p+({\frac{n-2}{p-2}}-\{\frac{n-2}{p-2}\} -1)(p-2)\\
& = &p+({\frac{n-p}{p-2}}-\{\frac{n-2}{p-2}\})(p-2)\\
& = &n-(p-2)\{\frac{n-2}{p-2}\}\\
& < &n.
\end{eqnarray*}
Therefore the tree $T$ does not touch all the faces of $M$. Hence $T$ is not an admissible proper tree. Similarly, when $m\geq\lceil\frac{n-2}{p-2}\rceil$, the number of faces in $\widetilde{F} > n$. But total number of faces in $M$ are $n$. This can not  happen. Hence $T$ is not an admissible proper tree. This proves the Lemma.

\smallskip

\noindent{\sc Proof of Lemma } \ref{lemma1}: First we prove, if $\frac{n - 2}{p - 2} = m$ is an integer then  $M$ has a proper tree on $m$ vertices. Let $K$ and $M$ be as in the statement of Lemma. We construct a vertex set $V(T)$ and an edge set $E(T)$ of a tree $T$ in $M$. For this, choose a vertex $v_1$ in $M$ and form $V(T) = \{v_1\}$. There is a facet $F_{v_1}$ in $K$ corresponding to $v_1$. Define a set $W = \{u \in V(K) \colon u$ is incident with $F_{v_1}\}$. At the $i^{th}$ step of construction, $1 < i \leq m$, form $V(T) = V(T) \cup {v_i}$ by choosing a vertex $v_i$ in $V(M)\setminus V(T)$ such that $v_{i}v_{j}$ is an edge in $M$ iff $j = i - 1$. Define $E(T) = \{v_{i-1}v_{i}\colon 1< i \leq m\}$. We get facets $F_{v_i}$ in $K$ corresponding to $v_i$ and put all the $p-2$ vertices of $F_{v_i}$ into $W$. In this construction there is no subset $U$ of $V(T)$ for which $U$ is equal to $V(F_{j}^*)$ for any facet $F_{j}^*$ in $M$. Thus, at the $m$-th step the graph $T := (V(T),E(T))$, which is a tree by construction, satisfies the conditions two and three in the Definition 1. Hence it is a proper tree. The number of elements in $W$ will be  $p+\underbrace{(p-2)+(p-2)+....+(p-2)}$ (repetitions of $(p-2) = (m-1))= p+(p-2)(m-1)= p + n - p = n$, since $m =\frac{n - 2}{p - 2} $. Hence by  Lemma \ref{lem2} the dual corresponding to $T$ bounds a $2-disk$ $D$ and $bd(D)$ is a Hamiltonian cycle in $K$. So, $T$ is an admissible proper tree. This proves the Lemma.

\section{The Steps for Searching a Separating Hamiltonian cycle in equivelar maps\,: }\label{algorithm}
The following steps may be implemented as a computer program to located separating Hamiltonian cycles\,:
\medskip
\begin{algo}\label{algo1}  Let $EG(K)$ be the edge graph of a $\{p,q\}$
equivelar map $K$ on a surface $S$, $\#V(EG(K)) = n$ and $w=\frac{n - 2}{p -
2}$ be an integer (by lemma \ref{lemm1} and  $lemma$ \ref{lem3}). Let $M$
be the set of all  $p$-gonal facet and $i$ denote the no of steps. We construct two set $D$ and $V$ as
follows. Choose an element $P_{0}\in M$. Define $D:=\{P_{0}\}$, $V :=
\{V(P_{0})\}$ and $i = 1$.

\begin{enumerate}

\item At next step, if $w = 1$ then, $n = p$ i.e. surface is a
$2-disk$ bounded by a $p$-gon. In this case, the $p$-gon itself a Hamiltonian cycle.
we stop hare.

\item Suppose $w > 1$. This follows $n > p$  as $n = p+(w-1)(p-2)$ and $p\geq 3$.
This follows $\#V < \#V(EG(K))$ i.e. there is a vertex which is not in $V$.
Then at next step, we go to the only one of the following steps.

\begin{enumerate}

\item Suppose there is a vertex $v\in V(EG(K))\setminus V$ and $v\in
V(P)$ with $V(P)\bigcap V=\{v_{1},v_{2}\}$ such that $E(P)\bigcap
E(P\_1)=\{\{v_{1},v_{2}\}\}$ where $P\_1\in D$. Then we take $D=D\bigcup \{P\}$,
$V = V\bigcup V(P)$ and we increase $i$ by 1 and go to the next step.

\item Suppose there is a vertex $v\in V(EG(K))\setminus V$ and $v\in V(P)$
with $V(P)\bigcap V=\emptyset$. Then we get a sequence of facets in order
$P_{1},P_{2},.....,P_{r}$ with the following properties-

\begin{enumerate}
\item $P_{i}\bigcap P_{i+1}$ is an edge for $1\leq i \leq r-1$
\item $P_{1}=P$
\item There exist only one facet $P\_1$ $(P\_1\in D)$ with $V(P\_1)\bigcap V(P_{r})=\{u_{1}, w_{1}\}$, where
$E(P\_1)\bigcap E(P_{r})=\{\{u_{1}, w_{1}\}\}$
\item For all $P_{i}'$s , $V(P_{i})\bigcap
V = \emptyset, i = 1, 2,\dots, r-1$.
\end{enumerate}
Let $W$ be the dual of $K$. Also $G$ denote the edge-graph of the dual map of
$(V,E(D),D)$ and $u$ denote the dual vertex of $P_{1}$ in $W$.\\
Then the above order sequence $P_{1},P_{2},.....,P_{r}$ of facets exist because-

\begin{enumerate}

\item There always exist a path $Q(v\rightarrow u)$ in $W$ with the following properties-

\begin{enumerate}

\item $v \in V(G)$
\item $u\in K$
\item $V(G)\bigcap V(Q)=\{v\}$

\end{enumerate}

as $G$ and $W$ are connected.

\item Dual of $Q$ in $K$ is $\{P_{1},P_{2},.....,P_{r}\}$

\end{enumerate}

Here we choose $P_{r}$ and we take $D = D \bigcup \{P_{r}\}$ and $V = V\bigcup V(P_{r})$.
Now we increase $i$ by 1 and go to the next step.

Hence after a step we will go to one of the above two steps until we get
a condition $V = V(EG(K))$ and $i = w$.
\end{enumerate}

 At last step, $\# D = w $ as $V = V(EG(K))$. Let $P_{1},P_{2},......,P_{w}$ be the facets in $D$ then
$P_{1}\bigcup P_{2}\bigcup, \dots, \bigcup P_{w}$ , $P_{j}\in D$ for
$1\leq j \leq w$ is a $2-disk$ and $\partial(\bigcup P_{i})$ is a
Hamiltonian cycle. We stop here.
\end{enumerate}

\end{algo}

\section{Proper tree in polyhedral maps\,:}\label{pmap}
\begin{lemma}\label{lem4} Let $T$ be a proper tree in a general polyhedral
map $M$ on a surface $S$. Then $T\bigcap F\neq \emptyset$ for any facet $F$
of $M$.
\end{lemma}

\noindent{\sc Proof of Lemma}\ref{lem4} Let $V(T)$ and $E(T)$ respectively
be the set of vertices and edges of the tree $T$. We construct two sets
$E$ and $\widetilde{F}$ as follows-

\begin{enumerate}

\item Choose a vertex $v_{1}\in V(T)$ of degree
$m$ and let $\{F_{1,1},F_{1,2},.....,F_{1,m}\}$ be the set of facets of $M$ such that
the facets $F_{1,1},F_{1,2},.....,F_{1,m}$ are adjacent to the vertex
$v_{1}$. Put $\widetilde{F}_1 = \{F_{1,1},F_{1,2},.....,F_{1,m}\}$ and $E_{1} = \{v_{1}\}$.

\item At 2nd step, choose a vertex $v_{2}$ other that $v_{1}$ where $\{v_{1},v_{2}\}$ is an
edge of $T$. Suppose the degree of $v_{2}$ is $l$. Then there are $l$ facets $F_{2,1},F_{2,2},.....,F_{2,l}$ adjacent to the vertex $v_{2}$. And here exactly two facets adjacent to $v_{1}$ and $v_{2}$ as $\{v_{1},v_{2}\}$ is an edge in polyhedral map. Hence, there are exactly $l-2$ no of facets of $M$ adjacent to $v_{2}$ which do not belong to the set $\widetilde{F}$. Hence, put all the new facets at $v_{2}$ in $\widetilde{F}_2$ and  $E_2 := E_1 \cup \{v_2\}$.

\item At a general step, say at $i^{th}$ step, choose a vertex $v_{i}\in V(T)\setminus E_{i -
1}$ such that for some $w\in E$, $\{w,v_{i}\}$ is an
edge of $T$ and assume $\deg(v_i)= t$. We define $E_i := E_{i-1}\cup
\{v_i\}$.

We claim, since $v_i \not\in E_{i - 1}$, there are exactly $r = (t-2)$ distinct facets
$F_{i,1},F_{i,2},\dots,F_{i,r}$ incident at $v_{i}$ and each are different from all
the facets of $\widetilde{F}$.
Suppose $r\neq(t-2)$. Then there are following possibilities:
\end{enumerate}
Suppose $r < t-2$ i.e. $t-r \geq 3$ then there are at least three facets
$F^{i}_{1},F^{i}_{2}$ and $F^{i}_{3}$ incident with $v_{i}$ and contained in $\widetilde{F}$.
Suppose $\{w,v_{i}\}$ is a common edge of
$F^{i}_{1}$ and $F^{i}_{2}$ and there exist a vertex $u\in E(T)$ such that
$u\in V(F^{i}_{3})$. This implies two distinct sub paths in $T$ have $u$ and $v_i$ as their end
vertices. These two paths would hence constitute a cycle in $T$ contradicting that $T$ is a tree. Hence $r \not < (t-2)$.
Arguing in similar way, we see that $r  > (t - 2)$ is also not possible. Hence $r = t-2$ i.e.
at $i^{th}$ step exactly $\deg(v_{i}) - 2$ new facets get added to the set $\tilde{F}$.  After $k^{th}$ step, the number of facets in
$\tilde{F} $ is $\deg(v_{1})+ \displaystyle\sum_{i = 2}^{k}(\deg(v_{i}) - 2)$, where number of elements in $V(T) = k$. In other words
number of elements in $\tilde{F}$ is  $\sum_{i = 1}^{ k}\deg(v_{i})-\sum_{i = 2}^{ k }2 = n + 2(k-1) - 2(k-1) = n$. So, after
$V(T) = k^{th}$  step $\tilde{F}$ will contain all the facets of the polyhedral map. Hence the tree $T$ touches all the facets of $M$.
This proves the Lemma.
\smallskip

\begin{lemma}\label{lem5} Let $K$ be a $n$ vertex polyhedral map and $M$
denote the dual polyhedron corresponding to $K$. Let $T$ be a
$k$ vertex proper tree in $M$. If $D$ denotes the subcomplex of $K$
which is dual of $T$ then $D$ is a $2$-disk and the boundary $\partial{D}$ of $D$ is a
Hamiltonian cycle in $EG(K)$.
\end{lemma}

\noindent{\sc Proof of Lemma } \ref{lem5}: Since $T$ has $k$ vertices, $D$
consists of $k$ facets $F_{1},F_{2},...,F_{k}$. Two facets in $D$ have an edge in common if the
corresponding vertices are adjacent in $T$. Here the set $D$ is a
collapsible $simplicial$ complex and hence it is a  $2-disk$. Since $T$ has
vertex of degree one, $bd(D)\neq\emptyset$, and being boundary complex
of a $2-disk$, it is a connected cycle. Observe that the number of edges
in $k~polygons$ are $\sum_{i = 1}^{ k }length(F_{i})$ and for each edge of
$T$ exactly $2~edges$ are identified. Hence the total number of edges which
remains unidentified in $D$ is $\sum_{i = 1}^{ k }degree(v_{i})-\sum_{i =
2}^{ k }2= n+2(k-1)-2(k-1)= n$. Hence the number of vertices in
$bd(D):=\partial D = n$. If the vertices $v_{1},v_{2}\in\partial D$ such
that $v_{1},v_{2}$ lie on a path of length $< n$ in $\partial D$ and $v_{1}=v_{2}$. This
means the facets $F_{1}$ and $F_{2}$ in $D$ with $v_{1}\in F_{1},v_{2}\in
F_{2},F_{1}\neq F_{2}$ and $F_{1}$ not adjacent to $F_{2}$. Thus there
exist a facet $F'$ in $D$ such that the vertex $u_{F'}$ in $T$
corresponding to $F'$ does not belong to the facet $F(v_{1})$ corresponding
to vertex $v_{1}$. But this contradicts that $T$ is a proper tree. Thus a
the cycle $\partial D$ contains exactly $n$ distinct vertices. Since
$\#V(K)=n$, $\partial D$ is a Hamiltonian cycle in $K$. This proves the Lemma.

\smallskip

\noindent {\sc Proof of Theorem } \ref{thm2}: The above lemma shows the if part. Conversely, let $K$ denote a polyhedral map and $H:=(v_{1},v_{2},....,v_{n}$) denote a contractible Hamiltonian cycle in $EG(K)$. Let $F_{1},F_{2},.....,F_{m}$ denote the facets of length $l_{i}=length(F_{i})$ such that $H = \partial (\bigcup_{j = 1}^m F_j = D)$. We claim that all the facets $F_i$ have their vertices on $H$. For, otherwise,  there will be identifications on the surface due to the hypothesis that $H$ is Hamiltonian. Thus, if $x$ denotes the number of facets in the  subpolyhedra $D$ which is topologically a 2-disc then the Euler characteristic relation gives us $1 = n-(\frac{\sum_{i = 1}^{ m }l_{i}-n}{2} + n) + m$ i.e. $\sum_{i = 1}^{ m }l_{i} = n + 2(m-1)$. In the edge graph of dual map $M$ of $K$, consider the graph corresponding to $D$ with $m$ vertices $u_{1}, u_{2},..., u_{m}$. This graph is a tree which is also a proper tree. This is so because $\sum_{i = 1}^{ m }length(F_{i}) = n + 2(m-1)$, where $l_{i}=length(F_{i})$ and there does not exist any subset $S_{1}$ of $ \{F_{1},F_{2},\ldots,F_{m}\}$ such that union of elements of $S_{1}$ is a $2$-disc subpolyhedra of $K$, whose boundary is a link of a vertex in $K$.
\hfill$\Box$

\begin{cor}\label{cor1}The edge graph $EG(K)$ of a $\{p, q\}$ equivelar map $K$ on a surface has a contractible Hamiltonian cycle if and only if the edge graph of corresponding dual map of $K$ has a proper tree.
\end{cor}

\noindent {\sc Proof of Corollary } \ref{cor1}: In the proof of above theorem \ref{thm2} we choose length $l_{i}= length(F_{i}) = p$ and $H:=(v_{1},v_{2},....,v_{n}$), where $H$ denotes a contractible Hamiltonian cycle in $EG(K)$. If $x$ denotes the number of $p$-gons in the disk resulting as a union of facets corresponding to $v_i$s then the Euler characteristic relation gives us $1 = n-(\frac{p\times x-n}{2} + n) + x$. Thus $x = \frac{n-2}{p-2}$. So that $m = \frac{n-2}{p-2}$. Hence, the edge graph $EG(K)$ of a $\{p, q\}$ equivelar map $K$ on a surface has a contractible Hamiltonian cycle. This proves the corollary.

\bigskip

\noindent {\sc Proof of Theorem } \ref{thm1}: The proof follows by corollary \ref{cor1}.
\hfill$\Box$

\end{document}